\documentclass{article}
\usepackage{a4wide}
\usepackage[centertags]{amsmath}
\usepackage{amsthm,amssymb,latexsym}
\usepackage{graphicx}
\usepackage{psfrag}

\newtheorem{thm}{Theorem}
\newtheorem{lem}{Lemma}
\newtheorem{obs}{Observation}  
\newtheorem{prop}{Proposition} 
\theoremstyle{definition}
\newtheorem{rem}{Remark}

\newcommand{\F}{{\cal F}}

\renewcommand{\S}{{\cal S}} 
\newcommand{\R}{{\mathbb R}}
\newcommand{\N}{{\mathbb N}}
\newcommand{\e}{\varepsilon}
\newcommand{\Rg}{{\R_{>0}}}

\newcommand{\sqsubsetdot}{\sqsubset\!\!\!\!\cdot\;}
\newcommand{\dcup}{\cup\!\!\!\!^\cdot\;\,}
\DeclareMathOperator{\conv}{conv}

\DeclareMathOperator{\gap}{gap}

\newenvironment{numenum}{%
\begin{enumerate}}{%
\end{enumerate}}

\newenvironment{alphenum}{%
\begin{enumerate}}{%
\end{enumerate}}

\begin{document}
\title{Kalai's squeezed 3-spheres are polytopal}
\author{\sc Julian Pfeifle\footnote{
Graduate student at the European graduate school
    `Combinatorics, Geometry, and Computation', 
supported by the Deutsche Forschungsgemeinschaft,
    grant GRK 588/1}\\
Dept.~of Mathematics, MA 6-2\\
 TU Berlin, D-10623 Berlin, Germany\\
{\tt pfeifle@math.tu-berlin.de}}
\date{October 19, 2001} 
\maketitle

\begin{abstract} 
In 1988, {\sc Kalai} \cite{Kalai88} extended a 
construction of {\sc Billera} and {\sc Lee}  
to produce {\it many} triangulated $(d-1)$-spheres. 
In fact, in view of upper bounds on the number of
simplicial $d$-polytopes by {\sc Goodman} and {\sc Pollack}
\cite{Good-Poll86,Good-Poll87}, he derived that 
for every dimension $d\ge 5$, {\em most} of these $(d-1)$-spheres
are not polytopal. However, for $d=4$, this reasoning fails.
We can now show that, as already conjectured by {\sc Kalai}, all of his
3-spheres are in fact polytopal. 

We also give a shorter proof for {\sc Hebble} and {\sc Lee}'s result \cite{Hebble-Lee00}
that the dual graphs of these $4$-polytopes are Hamiltonian.
\end{abstract}

\section{Introduction}
This paper is about triangulated spheres and the question whether or
not the members of a certain family of them are \emph{realizable},
that is, if they arise as boundary complexes of simplicial polytopes.
While for all two-dimensional spheres this is true by Steinitz'
Theorem, already one dimension higher there exist simplicial spheres
that cannot be realized in a convex way. The first example for
this---the so-called \emph{Br\"uckner sphere}---was found by {\sc
  Gr\"unbaum \& Sreedharan} in 1967, who realized that a certain
simplicial $3$-sphere on $8$ vertices does not represent the
combinatorial type of any $4$-polytope, contrary to what {\sc
  Br\"uckner} originally thought. (See \cite[Chapter 5]{Ziegler98} for a more
thorough discussion and references.)

In 1988, {\sc Kalai} extended a construction by {\sc Billera} and {\sc 
  Lee}, and showed that starting with $d=5$, there exist many more
simplicial $(d-1)$-spheres than simplicial $d$-polytopes, and that
therefore, in a very strong sense, \emph{most} simplicial spheres are
not realizable. In contrast, it is the main goal of this paper to show that 
all of {\sc Kalai}'s $3$-spheres \emph{do} arise as boundary complexes 
of simplicial $4$-polytopes.

In the remainder of this introduction, we present the context of these
constructions, including the known upper resp.~lower bounds for the
numbers of simplicial polytopes resp.~spheres.

The most important invariant of a $(d-1)$-dimensional simplicial
sphere $\S$ is its {\bf f-vector}
$f(\S)=(f_{-1},f_0,f_1,\dots,f_{d-1})$, where $f_i=f_i(\S)$ counts the
number of $i$-dimensional faces of $\S$, and $f_{-1}=1$.  In 1971,
{\sc McMullen} \cite{McMullen71} conjectured a characterization of the
$f$-vectors of boundary complexes of simplicial $d$-polytopes in terms
of an encoding of $f(\S)$, the so-called {\em $g$-vector}.  First
define the {\bf h-vector} $h(\S)=(h_0,h_1,\dots,h_d)$ of $\S$ by
\[
        h_k=\sum_{i=0}^{k} (-1)^{k-i}{d-i\choose d-k}
        f_{i-1},\qquad \text{for }k=0,1,\dots,d.
\]
The $h$-vector of any simplicial sphere
satisfies the {\bf Dehn-Sommerville equations} $h_k=h_{d-k}$ 
for $k=0,1,\dots,\lfloor d/2\rfloor$.
Now the {\bf g-vector} of $\S$ is
$g(\S)=(g_0,g_1,\dots,g_{\lfloor d/2\rfloor}),$ where
$g_0:=h_0=1$ and
\[
        g_k:=h_k-h_{k-1} \qquad\text{for }\quad
        k=1,2,\dots,\lfloor d/2\rfloor.
\]
We say that $g(\S)$ forms an {\bf M-sequence} if
$g_0=1$ and $g_{k-1}\ge\partial^k(g_k)$ for
$k=1,\dots,\lfloor d/2\rfloor$, where
\[
        \partial^k(g_k)={a_k-1\choose k-1} +
        {a_{k-1}-1\choose k-2}+\dots+
        {a_2-1\choose 1}+{a_1-1\choose 0},
\]
and the integers $a_k>a_{k-1}>\cdots>a_2>a_1\ge0$ are determined
by the {\bf binomial expansion}
\[
        g_k-1={a_k\choose k}+{a_{k-1}\choose k-1}+\dots
        +{a_2\choose 2}+{a_1\choose 1}
\]
of $g_k-1$ w.r.t. $k$.
See \cite[Chapter 8]{Ziegler98} for more details.
We can now state {\sc McMullen}'s conjecture:

\begin{thm} {\bf ($\mathbf g$-conjecture/theorem)}
An integer vector $g=(g_0,g_1,\dots,g_{\lfloor
d/2\rfloor})$ is the $g$-vector of the boundary complex of a simplicial
$d$-polytope $P$ if and only if it is an M-sequence.
\end{thm}

In the same year, 1979, {\sc Stanley} \cite{Stanley80} proved the  
necessity and {\sc Billera} and {\sc Lee} 
\cite{BL81} the sufficiency of \mbox{{\sc McMullen}'s} conditions.    
{\sc Stanley}'s proof that the $g$-vector of any simplicial   
polytope is an M-sequence used the Hard Lefschetz Theorem    
for the cohomology of projective toric varieties, but in 
the meantime a simpler proof by {\sc McMullen} using his  
{\em polytope algebra} is available. 
 
{\sc Billera} and {\sc Lee} invented an ingenious construction to 
produce, for every M-sequence $g$, a simplicial $d$-polytope with this
$g$-vector. Very briefly, they first find a shellable ball
$B$ as a collection of facets of a cyclic polytope $C$, 
such that the $g$-vector of $\partial B$ is the given M-sequence.
Then they construct a realization of $C$ and a
point $z$ that sees exactly the facets in $B$, and obtain
a realization of $\partial B$ as a simplicial polytope by taking
the vertex figure at~$z$ of $\conv(\{z\}\cup C)$. 

We next discuss {\sc Kalai}'s 1988 extension of their construction, by
which he built so many simplicial spheres that {\em most} of them (in
a sense to be made precise below) fail to be polytopal.  He achieved
this by giving a rule to produce many lists $I$ of $(d+1)$-tuples of
vertices, which span pure simplicial complexes $B(I)$. The underlying
space of every such complex turns out to be a simplicial, shellable
$d$-ball, which he called a {\bf squeezed ball}, and therefore the
boundary $S(I)$ of $B(I)$ is a simplicial $(d-1)$-sphere, a {\bf
  squeezed sphere}. {\sc Lee} shows in \cite{Lee00}
that {\sc Kalai}'s squeezed spheres are shellable.

Let $s(d,n)$ denote the number of
simplicial $(d-1)$-spheres, $sq(d,n)$ the number of
squeezed $(d-1)$-spheres, and $c(d,n)$ the number
of combinatorial types of simplicial $d$-polytopes with 
$n$ labeled vertices. {\sc Goodman}
and {\sc Pollack} \cite{Good-Poll86,Good-Poll87} derive the upper bound 
\begin{equation}\label{eq.polybound}
        \log c(d,n)\ \ \le\ \ d(d+1)n\log n  
\end{equation}
using  
a theorem of {\sc Milnor} that bounds the sum of the Betti numbers
of real algebraic varieties, 
while {\sc Kalai}'s squeezed spheres provide the
following lower bound for $s(d,n)$:
\begin{eqnarray*} 
        \log s(d,n)\;\ge\;\log sq(d,n) &\ge& \frac{1}{(n-d)(d+1)}
                {n-\lfloor(d+2)/2\rfloor \choose
                \lfloor(d+1)/2\rfloor } \\
              &=&\Omega(n^{\lfloor(d+1)/2\rfloor-1})
              \quad\text{for fixed }d.
\end{eqnarray*}
These bounds reveal that $\lim_{n\to\infty}c(d,n)/sq(d,n)=0$ for
$d\ge5$, which means that for $d\ge5$ most of {\sc Kalai}'s spheres
are not polytopal---there are simply too many of them. However, we
learn nothing for $d\le4$: We will prove in
Proposition~\ref{prop.count} below that $sq(4,n)\le 2^{n-5}n!$ for
$n\ge 5$, which is strictly less than the bound from
(\ref{eq.polybound}) for all $n\ge 5$.

\medskip The rest of the paper is organized as follows: In Section
\ref{sec.cyclic}, we collect some facts about \emph{cyclic polytopes},
an essential ingredient of our proof. In Section \ref{sec.real}, we
first present the details of {\sc Kalai}'s construction, and then show
how to realize any of his $3$-spheres as boundary complexes of
simplicial $4$-polytopes (Theorem \ref{thm.thm2}).
Finally, Section \ref{sec.hamil} uses the pictures constructed in
Section \ref{sec.real} to give a shorter proof of {\sc Hebble} and
{\sc Lee}'s result that the dual graphs of squeezed $3$-spheres are
Hamiltonian.

\section{Some facts on cyclic polytopes}\label{sec.cyclic}

The convex hull of $n$ distinct points on the {\bf moment curve}
$\mu_d:t\mapsto(t,t^2,\dots,t^d)$ in $\R^d$ is called a
$d$-dimensional {\bf cyclic polytope} with $n$ vertices. The
combinatorial type of this polytope is independent of the choice of
the $n$ points on the moment curve, and so one can talk about {\bf
  the} cyclic polytope $C_{d}(n)$. In fact, any $d$-dimensional order
$d$ curve also gives rise to the same combinatorial types of polytopes.

\medskip 
We switch from $d$ and $n$ to $d+1$ and $n+1$, and consider a
set $X=\{x_0=\mu(t_0),\dots,x_n=\mu(t_n)\}$ of $n+1$ distinct points
on the moment curve $\mu_{d+1}=:\mu$, ordered by their first
coordinates.  For any $f\subset\{0,1,\dots,n\}$, write $F_f$ for the
subset of $X$ indexed by~$f$, and $i(F)$ for the indices of a
subset~$F$ of $X$.  The supporting hyperplane~$H(F)$ of a
$(d+1)$-subset $F\subset X$ is given by
$H(F)=\{x\in\R^{d+1}:\gamma(F)\cdot x=-\gamma_0(F)\},$ where
$\gamma(F)=(\gamma_1(F),\dots, \gamma_{d+1}(F))\in\R^{d+1}$
and~$\gamma_0(F)\in\R$ are defined by
\begin{equation}\label{eq.prod}
        0\;=\;\prod_{i\in i(F)}(t-t_i) \;=\; \sum_{j=0}^{d+1}
        \gamma_j(F)t^j \;=\; \gamma_0(F)+\gamma(F)\cdot
        \mu(t).
\end{equation}
Observe that $\gamma_{d+1}(F)=1$;
we say that $\gamma(F)$ {\bf points upwards}.

\medskip
{\sc Gale}'s {\bf evenness criterion} tells us which $(d+1)$-subsets
$F$ of $X$ are vertex sets of facets of the cyclic polytope
$C=\conv(X)$: For any $i,j\in\{0,1,\dots,n\}\setminus i(F)$, the
number of elements of~$i(F)$ between~$i$ and~$j$ must be even.

\medskip
Define the {\bf end set} $W_{end}$ of $F_f\subset X$ to be the
right-most contiguous block $\{r_f+1,\dots,\max f\}$ of the indices
$f$ of $F$, where $r_f=\max\{i\in\N:i<\max f,\, i\notin f\}$.  Let $F$
be a facet of $C$ and take $x_j=\mu(t_j)\in X\setminus F$.  If the
cardinality of the end set of $F$ is odd, we get $\prod_{i\in
  i(F)}(t_j-t_i)<0$ because $j\notin i(F)$, and therefore
$\gamma(F)\cdot x_j< -\gamma_0(F)$.  Since $\gamma_{d+1}(F)=1$, we
conclude that the whole cyclic polytope $C$ is {\em below} $F$, and
call $F$ an {\bf upper facet} of $C$. If $\# W_{end}$ is even, we
analogously call~$F$ a {\bf lower facet} of $C$.  Finally, define an
{\bf outer normal vector} $\alpha(F)$ of any facet $F$ of $C$
by~$\alpha(F)=\gamma(F)$ resp.~$\alpha(F)=-\gamma(F)$ if $F$ is an
upper resp.~lower facet of $C$, and set~$\alpha_0(F)=-\gamma_0(F)$
resp.~$\alpha_0(F)=\gamma_0(F)$. By this, we obtain
$C\subset\{x\in\R^{d+1}:\alpha(F)\cdot x\le \alpha_0(F)\}$ for all
facets $F$ of $C$.

\section{Realizing Kalai's 3-spheres}\label{sec.real}

\subsection{{Kalai}'s idea}
First define a partial order $\preceq$ on ${\N\choose d+1}$ by
$\{i_1,i_2,\dots, i_{d+1}\}_< \preceq \{j_1,j_2,\dots, j_{d+1}\}_<$ if
$i_k\le j_k$ for every $k=1,\dots,d+1$. Here the notation
$A=\{a_1,\dots,a_r\}_<$ means that the elements of the set $A$ are
listed in increasing order. For the standard poset terminology used in
the following, see \cite{StanleyEnum1}. 

For an odd integer $d>0$ and $n\in\N$, let $\F_d(n)$ be the collection
of $(d+1)$-subsets of $[n]:=\{1,2,\dots,n\}$ of the form
$\{i_1,i_1+1\}\cup\{i_2,i_2+1\}\cup\cdots \cup\{i_e,i_e+1\}$, where
$e=(d+1)/2$, $i_1\ge1$, $i_e<n$, and $i_{j+1}\ge i_j+2$ for all
relevant $j$.  Let $I'$ be an initial set (order ideal) of $\F_d(n)$
with respect to the partial order $\preceq$ on ${\N\choose d+1}$.
Informally, $f'\preceq g'$ for $f',g'\in \F_d(n)$ if $f'$ arises from
$g'$ by pushing some elements in $g'$ to the left.

For even $d>0$, put $\F_d(n)=\{\{0\}\cup f':f'\in\F_{d-1}(n)\}
=:0*\F_{d-1}(n)$ with the induced partial order, and set $I:=0*I'$.

Finally, let $B(I)$ be the simplicial complex (the {\bf squeezed
  $d$-ball}) spanned by $I$, denote the boundary complex of $B(I)$ by
$S(I)$ (the {\bf squeezed $(d-1)$-sphere}), and do the same for $I'$.

\subsection{The structure of $3$-balls} \label{subsec.structure}

To specialize {\sc Kalai}'s construction to $d=4$, we first study
squeezed 3-balls. Take $n\ge4$ in $\N$, write $(i,j)$ for an element
$\{i,i+1,j,j+1\}\subset[n]$ of~$\F_3(n)$, and define the {\bf gap} of
$(i,j)\in\F_3(n)$ to be the number $j-i-2$ of integers between $i+1$
and $j$. From the fact that any two elements of $\F_3(n)$ with the
same gap are translates of each other and therefore
$\preceq$-comparable, we conclude that any~$\preceq$-antichain
in~$\F_3(n)$ can be linearly ordered by increasing gap, and denote
this order by~$\sqsubset$.  We remark that the difference between the
gaps of any two elements in a $\preceq$-antichain must be at least~2,
as otherwise the two elements would be $\preceq$-comparable. In
particular, the maximal number of elements of a $\preceq$-antichain in
$\F_3(n)$ is $\lceil (n-3)/2\rceil$.

\medskip   
Any order ideal $I'\subset\F_3(n)$ for $n\in\N$ is generated by the
set $G'=\{g_1', g_2',\dots, g_r'\}_\sqsubset$ of its maximal elements,
for some $r\le\lceil (n-3)/2\rceil$. By our discussion, the
$g_k'=(i_k,j_k)$ satisfy
\begin{numenum}  
\item $j_k\ge i_k+2$\, for $k=1,\dots,r,$ \quad and
\item $i_k>i_{k+1}$ and $j_k<j_{k+1}$\, for $k=1,\dots,r-1$.
\end{numenum} 

As an example, let $I'$ be the ideal generated by
$G'=\{(9,11),(8,12),(5,14),(2,17)\}_\sqsubset$:
\[  
\begin{picture}(0,0)%
\includegraphics{gens.pstex}%
\end{picture}%
\setlength{\unitlength}{3947sp}%
\begingroup\makeatletter\ifx\SetFigFont\undefined
\def\x#1#2#3#4#5#6#7\relax{\def\x{#1#2#3#4#5#6}}%
\expandafter\x\fmtname xxxxxx\relax \def\y{splain}%
\ifx\x\y   
\gdef\SetFigFont#1#2#3{%
  \ifnum #1<17\tiny\else \ifnum #1<20\small\else
  \ifnum #1<24\normalsize\else \ifnum #1<29\large\else
  \ifnum #1<34\Large\else \ifnum #1<41\LARGE\else
     \huge\fi\fi\fi\fi\fi\fi
  \csname #3\endcsname}%
\else
\gdef\SetFigFont#1#2#3{\begingroup
  \count@#1\relax \ifnum 25<\count@\count@25\fi
  \def\x{\endgroup\@setsize\SetFigFont{#2pt}}%
  \expandafter\x
    \csname \romannumeral\the\count@ pt\expandafter\endcsname
    \csname @\romannumeral\the\count@ pt\endcsname
  \csname #3\endcsname}%
\fi
\fi\endgroup
\begin{picture}(4815,552)(1306,-2491)
\put(1306,-2491){\makebox(0,0)[lb]{\smash{\SetFigFont{10}{12.0}{rm}1}}}
\put(4377,-2491){\makebox(0,0)[lb]{\smash{\SetFigFont{10}{12.0}{rm}12}}}
\put(4100,-2488){\makebox(0,0)[lb]{\smash{\SetFigFont{10}{12.0}{rm}11}}}
\put(3803,-2488){\makebox(0,0)[lb]{\smash{\SetFigFont{10}{12.0}{rm}10}}}
\put(3569,-2488){\makebox(0,0)[lb]{\smash{\SetFigFont{10}{12.0}{rm}9}}}
\put(3317,-2491){\makebox(0,0)[lb]{\smash{\SetFigFont{10}{12.0}{rm}8}}}
\put(1601,-2491){\makebox(0,0)[lb]{\smash{\SetFigFont{10}{12.0}{rm}2}}}
\put(1888,-2488){\makebox(0,0)[lb]{\smash{\SetFigFont{10}{12.0}{rm}3}}}
\put(5229,-2491){\makebox(0,0)[lb]{\smash{\SetFigFont{10}{12.0}{rm}15}}}
\put(4951,-2491){\makebox(0,0)[lb]{\smash{\SetFigFont{10}{12.0}{rm}14}}}
\put(5806,-2491){\makebox(0,0)[lb]{\smash{\SetFigFont{10}{12.0}{rm}17}}}
\put(6121,-2491){\makebox(0,0)[lb]{\smash{\SetFigFont{10}{12.0}{rm}18}}}
\put(2436,-2491){\makebox(0,0)[lb]{\smash{\SetFigFont{10}{12.0}{rm}5}}}
\put(2738,-2491){\makebox(0,0)[lb]{\smash{\SetFigFont{10}{12.0}{rm}6}}}
\put(4664,-2491){\makebox(0,0)[lb]{\smash{\SetFigFont{10}{12.0}{rm}13}}}
\end{picture}

\]     
Note that if $g'\sqsubset h'\in G'$, then~$g'$ is nested inside~$h'$
(possibly with overlap). From Figure \ref{fig.ideal} below, we will
read off the structure of the 3-ball $B(I')$ generated by $G'$, and
its boundary $S(I')$.

\begin{figure}[htp]
\[          
\input{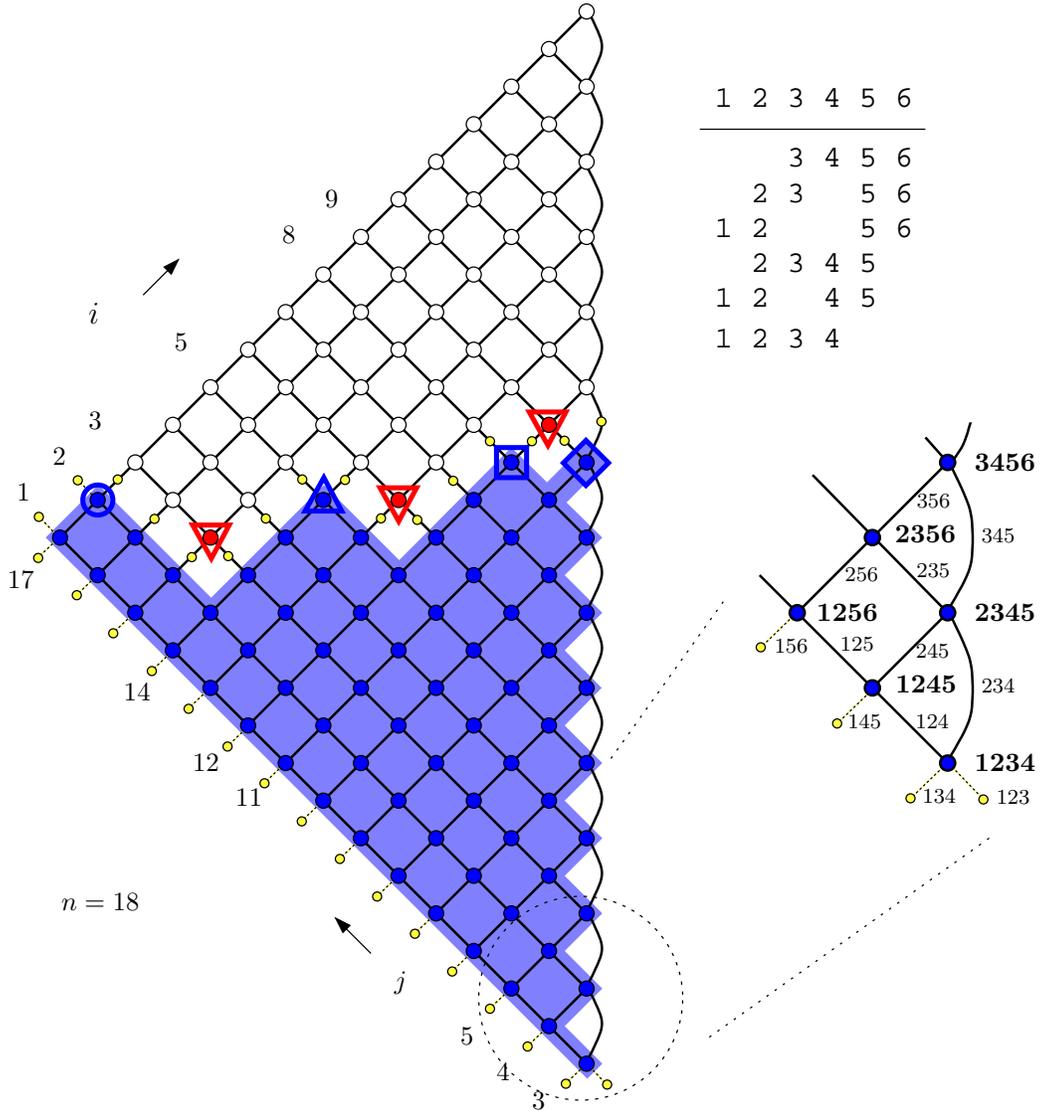}
\]        
\caption{The Kalai poset $\F_3(18)$. The shaded circles are
  the facets of the 3-ball $B(I')$ with
  generators~$G'=\{(9,11),(8,12),(5,14),(2,17)\}_\sqsubset$.  The
  minimal elements $H'$ of $\F_3(18)\setminus I'$ are marked by
  $\nabla$'s, and $\sqsubset$ orders the elements of $E'=G'\cup H'$
  from left to right (see Observation~\ref{obs.obs0}). Straight lines
  between facets correspond to $\prec$-covering relations between
  elements of $I'$, and straight and curved lines together to inner
  ridges of $B(I')$. The small circles are the facets of
  $S(I')=\partial B(I')$. The set of facets of the Kalai sphere $S(I)$
  is the union of $B(I')$ and $0*S(I')$.  }
\label{fig.ideal}
\end{figure}
 
Now put $\F_4(n)=0*\F_3(n)$ with the induced partial order, and $I=0*
I'$.  The 4-ball~$B(I)$ spanned by $I$ is a cone over the 3-ball
$B(I')$, whose boundary complex is the squeezed 3-sphere~$S(I)$.

\begin{prop}\label{prop.count}
There are at most  $2^{n-4}(n+1)!$ squeezed $3$-spheres with $n+1\ge 5$ 
labeled vertices. In particular, $\log sq(4,n)=\Theta(n\log n)$.
\end{prop}

\noindent{\em Proof.} By \cite[Prop.~3.3]{Kalai88}, distinct  
4-balls $B(I)$ whose vertices are labeled according to their
construction give rise to distinct 3-spheres $S(I)$ labeled in this
way, and distinct initial sets $I\subset\F_4(n)$ obviously induce
distinct such 4-balls. Every initial set $I$ is of the form $0*I'$ for
a unique order ideal~$I'\subset\F_3(n)$. Therefore, by relabeling
vertices, $sq(4,n+1)$ is at most $(n+1)!$ times the number of distinct
order ideals in $\F_3(n)$, depending on the combinatorial symmetries
of $S(I)$. By Figure \ref{fig.ideal}, every such order ideal can be
represented by a lattice path of length~$n-4$ taking steps only in the
positive~$i$- or negative~$j$-directions, and starting at
$(i,j)=(1,n-1)$. There are $2^{n-4}$ of these, and they all give rise
to distinct ideals.\hfill$\Box$

\subsection{A bird's-eye view of the realization construction}

Observe that by {\sc Gale}'s Evenness Criterion, every $f\in I$
corresponds to a lower facet~$F_f$ of a cyclic polytope.  By adapting
the ideas of {\sc Billera} and {\sc Lee}, we will now realize any
$S(I)$ as the boundary complex of a 4-polytope $P$ by appropriately
realizing a cyclic 5-polytope $C$, and choosing a viewpoint~$v$ close
to the negative $e_5$-axis that sees exactly the facets of $C$ in
$B(I)$.  The convex 4-polytope $P$ is then the vertex figure at $v$ of
$\conv(C\cup\{v\})$, and $S(I)$ its boundary.

\medskip
Specifically, let $\mu=\mu_5:\R\to\R^5$, $t\mapsto(t,t^2,\dots,t^5)$
be the moment curve in dimension 5.  Given an order ideal $I=0*I'$ in
$\F_4(n)$ where $n=\max\,\bigcup I$, we will execute the following steps:

\begin{enumerate}
\item Choose $N'>0$ and place $0=t_0<t_1<\cdots<t_n\in\R_{\ge0}$ such that 
\begin{equation}
\parbox{10cm}{
\begin{align*}
  \prod_{i\in f\setminus\{0\}}t_i<N' &\quad  \text{for all }
  f\in I,\qquad \text{and}\\
  \prod_{i\in f\setminus\{0\}}t_i>N' &\quad  \text{for all }
  f\in\F_4(n)\setminus I
\end{align*}
}\tag{S1}\label{eq.s1}
\end{equation}
Solutions for \eqref{eq.s1} exist with $t_1>0$ arbitrarily small. We
will find a solution for this system of inequalities by processing the
elements of $E'=G'\cup H'$ in $\sqsubset$-order, where $G'$ is the set
of $\preceq$-maximal elements of $I'$, and $H'$ is the set of
$\preceq$-minimal elements of $\F_3(n)\setminus I'$.
\item Make sure that the viewpoint to be defined will not see any
  upper facets of
  $C=C_5(n+1)=\conv\{0,\mu(t_1),\mu(t_2),\dots,\mu(t_n)\}$ that
  contain 0, by choosing $t_1>0$ so small that
\begin{equation}
   t_1t_{n-2}t_{n-1}t_n<N'.\tag{S2}\label{eq.inequpper}
\end{equation}
\item Choose $\e$, with $0<\e<t_1$, so small that for all $e,f\in\F_4(n)$,
\begin{equation}
   e\prec f \quad\Longrightarrow\quad \gamma(F_e)\cdot\mu(\e)
   \;<\;\gamma(F_f)\cdot\mu(\e).\tag{S3}\label{eq.s3}
\end{equation}
\item Choose $\e>0$ even smaller, if necessary, such that the
  viewpoint $v:=\mu(\varepsilon)-\e N'e_5$ satisfies
\begin{align}
  \alpha(F)\cdot v>\alpha_0(F) \quad& \text{for }f_F\in I,\notag \\
  \alpha(F)\cdot v<\alpha_0(F) \quad& \text{for all lower facets }F
  \text{ of } C \text{ such that }f_F\notin I,\tag{S4}\label{eq.s4} \\
  \alpha(F)\cdot v<\alpha_0(F) \quad& \text{for all upper facets }
  F \text{ of }C,\notag
\end{align}
where $\alpha(F)$ is the outer normal vector of $F$ we defined at the end
of Section \ref{sec.cyclic}.
\end{enumerate}

\noindent We conclude that $v$ sees exactly the facets of $C$  in 
$B(I)$, and obtain $S(I)$ as above.

\subsection{How to realize Kalai's 3-spheres}

We will now give the details of the construction and prove the following
theorem.

\begin{thm}\label{thm.thm2}
  Every squeezed $3$-sphere $S(I)$ given by an order ideal $I$ in the
  poset $(\F_4(n),\preceq)$  with $n\ge\max\,\bigcup I$ can be realized
  as the boundary complex of a simplicial, convex $4$-polytope.
\end{thm}

\begin{rem}
  The construction shows the stronger result that every squeezed
  $4$-ball $B(I)$ can be realized as a regular triangulation of a convex
  $4$-polytope.
\end{rem}

To prove Theorem \ref{thm.thm2}, given an ideal $I\subset\F_4(n)$, we
may assume that $n=\max\,\bigcup I$ since $\F_4(n)\subseteq\F_4(n')$
for $n\le n'$. By definition, every order ideal $I\subset\F_4(n)$ has
the form $I=0*I'$, where $I'=\langle G'\rangle\subset\F_3(n)$ is
generated by its maximal elements $G'=\{g_1',g_2',\dots,g_r'\}$ with
$g_k'=(i_k,j_k)$. Choose~$N'>0$, introduce $n$ variable points
$0<t_1<t_2<\dots<t_n$ in~$\Rg$, and consider the set~$H'$ of
$\preceq$-minimal elements of $\F_3(n)\setminus I'$.

\begin{obs}\label{obs.obs0}
  Consider any two consecutive elements $e'=(i,j)\sqsubsetdot
  f'=(k,\ell)$ of a $\sqsubset$-ordered
  \mbox{$\preceq$-antichain~$G'$} of~$\F_3(n)$.  Then the unique
  $\prec$-minimal element $m'$ in $\F_3(n)\setminus\langle G'\rangle$
  with $\gap(e')<\gap(m')<\gap(f')$ exists and is $m'=(k+1,j+1)$.  In
  particular, the number of $\prec$-minimal elements
  in~$\F_3(n)\setminus\langle G'\rangle$ is no greater than
  $\lfloor(n-3)/2\rfloor$.  \hfill$\Box$
\end{obs}

\noindent\emph{Sketch of proof.} The first statement follows by
inspection of Figure~\ref{fig.ideal}. For the second assertion, note
that the set $H'$ has maximal cardinality if
$G'=\{(i,n-i):i=1,2,\dots,\lceil(n-3)/2\rceil\}$.\hfill$\Box$

\medskip

Using Observation \ref{obs.obs0}, we linearly order $E'=G'\cup H'$ by
$\sqsubset$, see Figure \ref{fig.ideal}. To carry out Step~1 of our
program, first choose some small $\delta>0$. Our goal is to place the
$t$'s in $\Rg$ such that 
\begin{equation}\label{eq.s1prime}
   \prod_{i\in g'}t_i= N'-\delta \quad\text{for }g'\in G'
   \qquad\qquad\text{and}\qquad\qquad
   \prod_{i\in h'}t_i= N'+\delta \quad\text{for }h'\in H'.\tag{S$1'$}
\end{equation}

\begin{obs}\label{obs.card}
  The cardinality of $E'=G'\dcup H'$ is at most $n-3$. In particular,
  there are fewer equalities in \eqref{eq.s1prime} 
  than there are variables.
\end{obs}

\noindent {\em Proof.} Because $n=\max\,\bigcup I$, the largest element of
$(E',\sqsubset)$ is in $G'$. Using Observation \ref{obs.obs0} again,
\[
   \# E'=\# G' + \# H'\le\left\lceil\frac{n-3}{2}\right\rceil + 
   \left\lfloor\frac{n-3}{2}\right\rfloor = n-3,
\]
which proves Observation \ref{obs.card}. \hfill$\Box$

\medskip

We now begin the construction by placing the $t$'s corresponding to
the $\sqsubset$-smallest element of $E'$ in such a way in $\Rg$
that~\eqref{eq.s1prime} is satisfied. This is clearly possible. The
general step of constructing a solution to~\eqref{eq.s1prime} is based
on the following lemma.

\begin{lem}\label{lem.lemmaone}
Let $e'=(i,j)\sqsubsetdot f'=(k,\ell)$ be two consecutive elements of $E'$.
\begin{alphenum}
\item If $e'\in G'$ and $f'\in H'$, then $0<k\le i$ and $\ell=j+1$. If
  $e'\in H'$ and $f'\in G'$, then $k=i-1$ and $j\le\ell<n$. (See
  Figure~\ref{fig.ideal}.)
\item\label{lem.lemmaoneb} Suppose that the $\{t_i\}_{i\in e'}$ have
  been placed already, but not all $\{t_j\}_{j\in f'}$. Then these
  latter $t$'s may be placed in such a way in $\Rg$ that
  $0<t_k<t_{k+1}<t_\ell<t_{\ell+1}$, and the equality
\begin{equation}\label{eq.M}
   t_kt_{k+1}t_\ell t_{\ell+1} = M
\end{equation}
is satisfied, where $M:=N'-\delta$ if $f'\in G'$ and $M:=N'+\delta$ if
$f'\in H'$.
\end{alphenum}
\end{lem}

\noindent {\em Sketch of proof for} (\ref{lem.lemmaoneb}). 
Suppose that $e'\in G'$ and $f'\in H'$.  We then have the following
situation:
\[
\begin{picture}(0,0)%
\includegraphics{case1.pstex}%
\end{picture}%
\setlength{\unitlength}{3947sp}%
\begingroup\makeatletter\ifx\SetFigFont\undefined
\def\x#1#2#3#4#5#6#7\relax{\def\x{#1#2#3#4#5#6}}%
\expandafter\x\fmtname xxxxxx\relax \def\y{splain}%
\ifx\x\y   
\gdef\SetFigFont#1#2#3{%
  \ifnum #1<17\tiny\else \ifnum #1<20\small\else
  \ifnum #1<24\normalsize\else \ifnum #1<29\large\else
  \ifnum #1<34\Large\else \ifnum #1<41\LARGE\else
     \huge\fi\fi\fi\fi\fi\fi
  \csname #3\endcsname}%
\else
\gdef\SetFigFont#1#2#3{\begingroup
  \count@#1\relax \ifnum 25<\count@\count@25\fi
  \def\x{\endgroup\@setsize\SetFigFont{#2pt}}%
  \expandafter\x
    \csname \romannumeral\the\count@ pt\expandafter\endcsname
    \csname @\romannumeral\the\count@ pt\endcsname
  \csname #3\endcsname}%
\fi
\fi\endgroup
\begin{picture}(4614,915)(1339,-2536)
\put(1576,-2491){\makebox(0,0)[lb]{\smash{\SetFigFont{10}{12.0}{rm}$t_k$}}}
\put(3016,-2491){\makebox(0,0)[lb]{\smash{\SetFigFont{10}{12.0}{rm}$t_i$}}}
\put(1801,-2491){\makebox(0,0)[lb]{\smash{\SetFigFont{10}{12.0}{rm}$t_{k+1}$}}}
\put(3286,-2491){\makebox(0,0)[lb]{\smash{\SetFigFont{10}{12.0}{rm}$t_{i+1}$}}}
\put(4636,-2491){\makebox(0,0)[lb]{\smash{\SetFigFont{10}{12.0}{rm}$t_j$}}}
\put(4951,-2491){\makebox(0,0)[lb]{\smash{\SetFigFont{10}{12.0}{rm}$t_{j+1}=t_\ell$}}}
\put(5806,-2491){\makebox(0,0)[lb]{\smash{\SetFigFont{10}{12.0}{rm}$t_{\ell+1}$}}}
\put(5896,-1771){\makebox(0,0)[lb]{\smash{\SetFigFont{10}{12.0}{rm}$c$}}}
\put(5266,-1771){\makebox(0,0)[lb]{\smash{\SetFigFont{10}{12.0}{rm}$c_0$}}}
\put(2971,-1771){\makebox(0,0)[lb]{\smash{\SetFigFont{10}{12.0}{rm}$b_0$}}}
\put(1891,-1771){\makebox(0,0)[lb]{\smash{\SetFigFont{10}{12.0}{rm}$b$}}}
\put(1576,-1771){\makebox(0,0)[lb]{\smash{\SetFigFont{10}{12.0}{rm}$a$}}}
\end{picture}
 
\] 
It is straightforward to verify that
for any $0<k\le i$, the points
$a,b,c$ may be placed in such a way that $0<a<b<b_0<c_0<c$ and
$abc_0c=N'+\delta$. Similarly, if $e'\in H'$ and $f'\in G'$, 
\[
\begin{picture}(0,0)%
\includegraphics{case2.pstex}%
\end{picture}%
\setlength{\unitlength}{3947sp}%
\begingroup\makeatletter\ifx\SetFigFont\undefined
\def\x#1#2#3#4#5#6#7\relax{\def\x{#1#2#3#4#5#6}}%
\expandafter\x\fmtname xxxxxx\relax \def\y{splain}%
\ifx\x\y   
\gdef\SetFigFont#1#2#3{%
  \ifnum #1<17\tiny\else \ifnum #1<20\small\else
  \ifnum #1<24\normalsize\else \ifnum #1<29\large\else
  \ifnum #1<34\Large\else \ifnum #1<41\LARGE\else
     \huge\fi\fi\fi\fi\fi\fi
  \csname #3\endcsname}%
\else
\gdef\SetFigFont#1#2#3{\begingroup
  \count@#1\relax \ifnum 25<\count@\count@25\fi
  \def\x{\endgroup\@setsize\SetFigFont{#2pt}}%
  \expandafter\x
    \csname \romannumeral\the\count@ pt\expandafter\endcsname
    \csname @\romannumeral\the\count@ pt\endcsname
  \csname #3\endcsname}%
\fi
\fi\endgroup
\begin{picture}(4614,870)(1339,-2491)
\put(2341,-1771){\makebox(0,0)[lb]{\smash{\SetFigFont{10}{12.0}{rm}$a_0$}}}
\put(5581,-1771){\makebox(0,0)[lb]{\smash{\SetFigFont{10}{12.0}{rm}$b$}}}
\put(4636,-1771){\makebox(0,0)[lb]{\smash{\SetFigFont{10}{12.0}{rm}$b_0$}}}
\put(5581,-2446){\makebox(0,0)[lb]{\smash{\SetFigFont{10}{12.0}{rm}$t_\ell$}}}
\put(4321,-2446){\makebox(0,0)[lb]{\smash{\SetFigFont{10}{12.0}{rm}$t_j$}}}
\put(5896,-1771){\makebox(0,0)[lb]{\smash{\SetFigFont{10}{12.0}{rm}$c$}}}
\put(1576,-1771){\makebox(0,0)[lb]{\smash{\SetFigFont{10}{12.0}{rm}$a$}}}
\put(1621,-2446){\makebox(0,0)[lb]{\smash{\SetFigFont{10}{12.0}{rm}$t_k$}}}
\put(2026,-2446){\makebox(0,0)[lb]{\smash{\SetFigFont{10}{12.0}{rm}$t_{k+1}=t_i$}}}
\put(3016,-2446){\makebox(0,0)[lb]{\smash{\SetFigFont{10}{12.0}{rm}$t_{i+1}$}}}
\put(4546,-2446){\makebox(0,0)[lb]{\smash{\SetFigFont{10}{12.0}{rm}$t_{j+1}$}}}
\put(5851,-2446){\makebox(0,0)[lb]{\smash{\SetFigFont{10}{12.0}{rm}$t_{\ell+1}$}}}
\end{picture}
 
\] 
for any  $j\le\ell<n$ we may place $a,b,c$ such that $0<a<a_0<b_0<b<c$ and
$aa_0bc=N'-\delta$. \hfill$\Box$

\medskip 

We now complete Step 1 by applying Lemma \ref{lem.lemmaone} to all
members of $E'$ in $\sqsubset$-order. The definition of $\preceq$ tells
us that because the $f'\in E'$ satisfy~\eqref{eq.s1prime}, in fact all
$f\in\F_4(n)$ satisfy the system~\eqref{eq.s1}.

\medskip
If in Step 1 we encountered some $e'\in E'$ with $1\in e'$, then
necessarily $e'=\{1,2,n-1,n\}\in G'$, which imposed the inequality
$t_1t_2t_{n-1}t_n<N'$. This inequality in turn remains satisfied if we
choose $t_1$ even small enough to verify~(\ref{eq.inequpper}). If $1\notin
e'$ for all $e'\in E'$, we are free to do the same.  We have completed
Step 2, and place any remaining unassigned $t$'s such that
$0=t_0<t_1<\cdots<t_n$.

\begin{obs}\label{obs.obs2} \begin{alphenum}
\item $\gamma_0(F_f)=0$ for any $5$-element 
subset $f\subset \{0,1,\dots,n\}$ that contains $0$.
\item\label{obs.obs2b} For all choices of $t_1<\cdots<t_n$, 
  one can find $\e>0$ small enough such that the
  implication~\eqref{eq.s3} holds for all $f,g\in\F_4(n)$.
\end{alphenum}
\end{obs}

\noindent{\em Proof of} (\ref{obs.obs2b}). 
The definition (\ref{eq.prod}) of the $\gamma$'s implies that for
$f=\{0,s_1,\dots,s_4\}$,
\begin{equation}\label{eq.epsilon}
        \gamma(F_f)\cdot \mu(\e)=\e(\e-s_1)\cdots(\e-s_4)= 
        \e s_1s_2s_3s_4\pm o(\e). 
\end{equation}
This means that $\gamma(F_f)\cdot \mu(\e)<\gamma(F_g)\cdot \mu(\e)$ by
definition of $\prec$, for $\e$ small enough.\hfill$\Box$

\medskip Take $0<\e<t_1$ as in Observation \ref{obs.obs2}(\ref{obs.obs2b}),
tentatively set $z:=\mu(\e)$, and let $f\in\F_4(n)$. If~$f\in I$,
there exists some $g\in G:=0*G'$ with~$f\preceq g$, and by
\eqref{eq.epsilon}, we have
\[
  \gamma(F_f)\cdot z \;\le\; \gamma(F_g)\cdot z =\e\prod_{i\in
  g\setminus\{0\}}t_i+O(\e^2) = \e(N'-\delta)\pm o(\e).
\]
If $f\notin I$, then there is some $h\in H:=0*H'$ with $f\succeq h$,
and we obtain in a similar way that
\[
  \gamma(F_f)\cdot z \;\ge\; \e(N'+\delta)\pm o(\e).
\]
Thus, we finally choose $0<\e<t_1$ so small that with $z:=\mu(\e)$ and
$N:=\e N'$, we have $\gamma(F_f)\cdot z<N$ for $f\in I$, and $\gamma(F_f)
\cdot z>N$ for $f\notin I$. Step 3 is now complete.

\medskip

We proceed to verify that $v:=\mu(\e)-\e N'e_5=z-Ne_5$ satisfies the
inequalities \eqref{eq.s4}.  For this, recall that all $F_f$ with
$f\in \F_4(n)$ satisfy {\sc Gale}'s Evenness Criterion, which means
that~$\F_4(n)$ is exactly the set of lower facets of the cyclic
polytope $C=\conv(X)$ that contain $x_0=0$.  However, {\em
  any}~$F\subset X$ of odd cardinality satisfying {\sc Gale}'s
Evenness Criterion with even end-set must contain~$0$, and we conclude
that $\F_4(n)$ is in fact the set of {\em all} lower facets of $C$.

\medskip
Recall from Section \ref{sec.cyclic} that $\alpha(F)=\gamma(F)$ and
$\alpha_0(F)=-\gamma_0(F)$ if $F$ is an upper facet of $C$, and that
$\alpha(F)=-\gamma(F)$ and $\alpha_0(F)=\gamma_0(F)$ if $F$ is a lower
facet of $C$.  We  and discuss all facets $F_f$ of $C$ in turn:

\medskip
\noindent{\em Lower facets of $C$:}  
\begin{itemize}
\item If $f\in I\subset\F_4(n)$, then by 
construction $\gamma(F_f)\cdot z<N$, and this implies
$\gamma(F_f)\cdot v<0$ (remember that $\gamma_5(F)=1$ for all $F$) and
$\alpha(F_f)\cdot v>0=\alpha_0(F_f)$, which means that 
$F_f$ is visible from $v$.

\item If $f\in\F_4(n)\setminus I$,  we conclude from
$\gamma(F_f)\cdot z>N$ that $\alpha(F_f)\cdot
v<0=\alpha_0(F_f)$, which says that~$F_f$ is not visible
from $v$. 
\end{itemize}

\noindent{\em Upper facets of $C$:}
\begin{itemize}
\item If $0\not\in f=\{s_1,\dots,s_5\}$, then
(\ref{eq.prod}) and $\e<t_1$ imply $\gamma(F_f)\cdot
z+\gamma_0(F_f)=\prod_{i=1}^5(\e-s_i)<0$, and
\[
        \alpha(F_f)\cdot v = \gamma(F_f)\cdot v =
        \gamma(F_f)\cdot z-N 
        < -\gamma_0(F_f)-N < -\gamma_0(F_f) =
        \alpha_0(F_f).
\]
\item If $0\in f$, then $\gamma_0(F_f)=0$ and
  $f=\{0,1\}\cup\{i,i+1\}\cup\{n\}$ with $2\le i\le n-2$.  By
  inequality~(\ref{eq.inequpper}) and the definition of $\prec$, we
  conclude that necessarily $\gamma(F_f)\cdot z<N$ and
\[
        \alpha(F_f)\cdot v=
        \gamma(F_f)\cdot z-N < 0 = \alpha_0(F_f).
\]
\end{itemize}

\noindent We have verified the inequalities \eqref{eq.s4} and completed
the proof of Theorem \ref{thm.thm2}.\hfill $\Box$

\begin{rem}
  A referee has suggested to extend this construction to boundaries of
  more general even-dimensional squeezed balls. However, so far we
  have only been able to realize odd-dimensional squeezed spheres
  directly modeled on the $3$-dimensional ones, and leave this as an
  open problem.
\end{rem}

\section{A shorter proof that squeezed $3$-spheres are 
  Hamiltonian} \label{sec.hamil}

In 1973, {\sc Barnette}~\cite{Rosenfeld-Barnette73} conjectured that
all simple $4$-polytopes admit a Hamiltonian circuit.
In~\cite{Hebble-Lee00}, {\sc Hebble} and {\sc Lee} prove that squeezed
$3$-spheres are (dual) Hamiltonian by explicitly constructing a
Hamiltonian circuit in the dual graph; however, their proof goes
through extensive case analysis. A referee has suggested that it might
be possible to obtain a simpler proof of this result. In this section,
we follow his or her suggestion and obtain a ``proof by picture'' with fewer
case distinctions, which moreover only depend on parity conditions.

\begin{thm}\label{thm.hl} \emph{({\sc Hebble} and {\sc Lee}, 
    2000~\cite{Hebble-Lee00})} The dual graph of any Kalai
  $4$-polytope $S(I)$ admits a Hamiltonian circuit. In particular,
  the polars of these $4$-polytopes satisfy {\sc Barnette}'s conjecture.
\end{thm}

\noindent\emph{Proof.} Recall from Section~\ref{subsec.structure} that
the set of facets of $S(I)$ is $B(I')\cup (0*S(I'))$. We continue to
write $(i,j)=\{i,i+1,j,j+1\}$ for facets of $S(I)$ in $B(I')$, and
introduce the notation $(i+\frac12,j):=\{0,i+1,j,j+1\}$ and
$(i,j+\frac12):=\{0,i,i+1,j+1\}$ for facets of $S(I)$ in $0*S(I')$.
Also, recall from Section~\ref{subsec.structure} the definition of the
order relations $\preceq$ and $\sqsubset$, and number the set $G'$ of
$\preceq$-maximal elements $(i_k,j_k)$ of $B(I')$ in ascending
$\sqsubset$-order, starting with $k=1$.

We start our Hamiltonian circuit in the dual graph of $S(I)$ at the
facet $(i_0,j_0)=(1,3)=\{1,2,3,4\}\in B(I')$. While walking through
the other facets of $B(I')$, we will also pick up the facets of the
form $(i+\frac12,j)$ and $(i,j+\frac12)$ with $i,j\ge 1$ of $S(I')$,
and then return to $(1,3)$ via the set of facets $\{(0,j):2\le j\le
n-1\}$. We will also use the difference operators $\Delta
j_k=j_{k+1}-j_k$ and $\Delta i_k=i_{k+1}-i_k$. In our circuit, we
repeatedly go through certain steps, and in the figures we will mark
the end of one step and the beginning of the next by a square.  In all
steps, if all facets in $G'$ are processed, go to step \emph{Down} (and 
then to \emph{Finish}).

\begin{enumerate}
\item\emph{Over the top:} Start at $(i_0,j_0)=(1,3)$. If $j_1-j_0$ is
  odd, continue as in Figure~\ref{fig:over-down}(a).  If $j_1-j_0$ is
  even, proceed as in Figure~\ref{fig:over-down}(b). In both cases, go
  on until $(i_1+\frac12,j_1)$. Set $k=1$, and go to step \emph{Down}.
  \begin{figure}[htp] 
    \begin{center}
      \psfrag{i}{$i$} 
      \psfrag{j}{$j$}
      \psfrag{i1}{\small $i_1$}
      \psfrag{i2}{\small $i_2$}
      \psfrag{j1}{\small $j_1$}
      \psfrag{zero}{\small$(i_0,j_0)$}
      \psfrag{(a)}{\small (a) $j_1-j_0$ odd}
      \psfrag{(b)}{\small (b) $j_1-j_0$ even}
      \includegraphics[width=.75\textwidth]{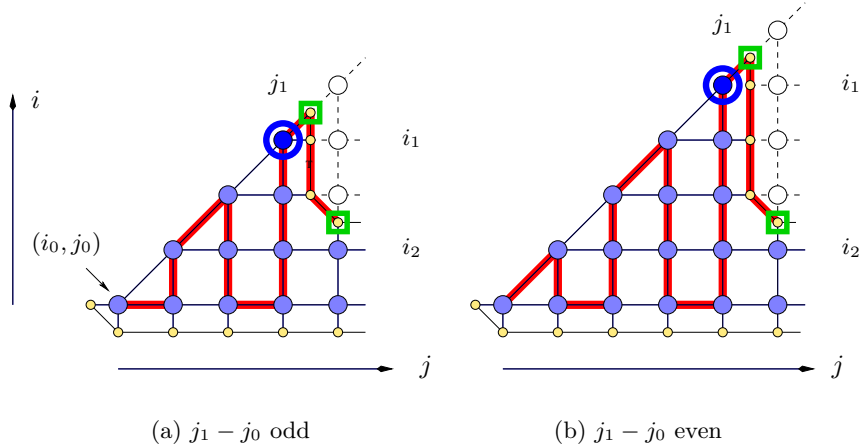} 
      \caption{Steps \emph{Over the top} and \emph{Down}. The circled facet is 
        $(i_1,j_1)$, the upper $\Box$ represents $(i_1+\frac12,j_1)$,
        and the lower $\Box$ is $(i_2+\frac12,j_1+1)$.}
      \label{fig:over-down}   
    \end{center}
  \end{figure}
  
\item\emph{Down:} If there are no more generators to be processed, go
  down along the facets $\{(i_\ell,j_k+\frac12):\ell=k,k-1,\dots,1\}$
  and continue with step \emph{Finish}. Otherwise, if $\Delta i_k>0$,
  continue downwards as in Figure~\ref{fig:over-down} until
  $(i_{k+1}+\frac12,j_k+1)$. If $i_{k+1}=i_k$, do nothing. In both
  cases, increment $k$ by $1$, and continue to step \emph{Across}.
  
\item\emph{Across:} If $\Delta j_k$ is even, continue as in
  Figure~\ref{fig:across-even}(a).  If $\Delta j_k$ is odd and not $1$
  and $i_{k+1}-i_0$ is even, continue as in
  Figure~\ref{fig:across-even}(b); if $\Delta j_k\ne1$ and $i_{k+1}-i_0$
  are both odd, as in Figure~\ref{fig:across-even}(c).

  \begin{figure}[htp]
    \begin{center}
      \psfrag{i0}{\small $i_0$}
      \psfrag{ik1}{\small $i_{k+1}$}
      \psfrag{jk}{\small $j_k$}
      \psfrag{jk1}{\small $j_{k+1}$}
      \psfrag{jeven}{\small (a) $\Delta j_k$ even}
      \psfrag{joddiodd}{\small (b) $\Delta j_k$ odd, $i_{k+1}-i_0$ odd}
      \psfrag{joddieven}{\small (c) $\Delta j_k$ odd, $i_{k+1}-i_0$ even}
      \includegraphics[width=.85\textwidth]{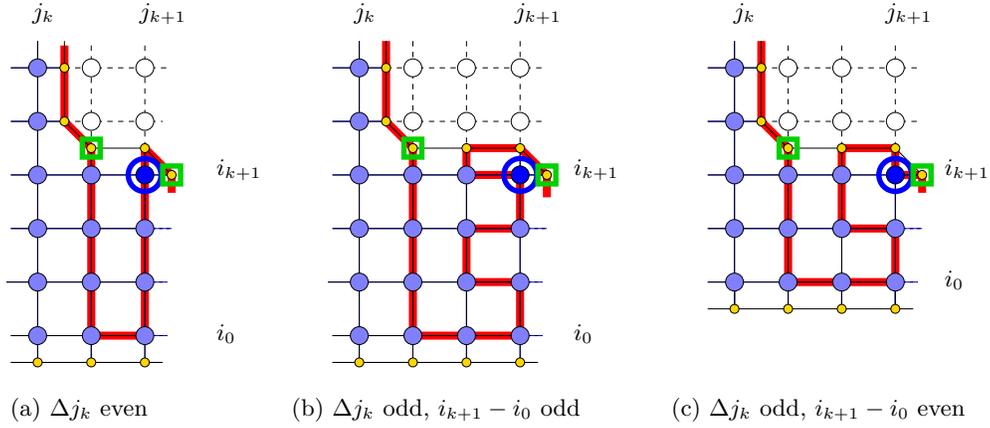}      
      \caption{Step \emph{Across} in case $\Delta j_k$ is even. The circled facet is
        $(i_{k+1},j_{k+1})$. }
      \label{fig:across-even}
    \end{center} 
  \end{figure}  
  
  If $\Delta j_k=1$ and $\Delta i_{k+1}$ is even, proceed as in
  Figure~\ref{fig:across-odd}(a), if $\Delta i_{k+1}$ is odd, as in
  Figure~\ref{fig:across-odd}(b).  In any case, increment $k$ by one,
  and repeat from step \emph{Down} or \emph{Across} as necessary,
  depending on whether the facet surrounded by a dashed circle in Figure
  \ref{fig:across-odd} is in $G$ or not.

  \begin{figure}[htp]
    \begin{center}
      \psfrag{ik1}{\small $i_{k+1}$}
      \psfrag{ik2}{\small $i_{k+2}$}
      \psfrag{jk}{\small $j_k$}
      \psfrag{jk1}{\small $j_{k+1}$}
      \psfrag{ieven}{\small (a) $\Delta i_{k+1}$ even}
      \psfrag{iodd}{\small (b) $\Delta i_{k+1}$ odd}
      \includegraphics[width=.6\textwidth]{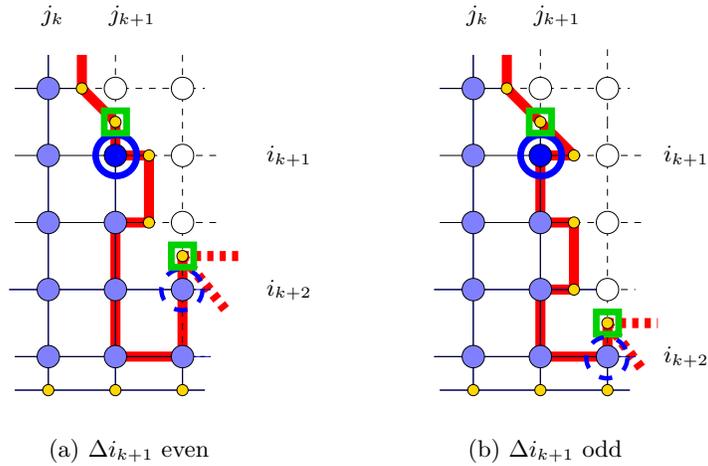}      
      \caption{Step \emph{Across} in case $\Delta j_k=1$. The circled facet is
        $(i_{k+1},j_{k+1})$. Depending on whether the facet surrounded
        by a dashed circle in Figure \ref{fig:across-odd} is in $G$ or
        not, the next step will be \emph{Down} or \emph{Across}, respectively.}
      \label{fig:across-odd}
    \end{center}
  \end{figure}  
 
\item\emph{Finish:} Now the only thing left to do is to return to
  $(1,3)$ via the set of facets $\{(0,j):n-1\ge j\ge 2\}$, as
  in Figure~\ref{fig:finish}.
 
  \begin{figure}[htp]
    \begin{center}
      \includegraphics[width=.5\textwidth]{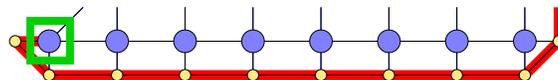}      
      \caption{Step \emph{Finish}.}
      \label{fig:finish}
    \end{center}
  \end{figure}  

\end{enumerate}
\noindent This completes the proof of Theorem \ref{thm.hl}.\hfill$\Box$

\section{Acknowledgements} 
It is a great pleasure to thank {\sc G\"unter M.~Ziegler} for
suggesting this problem, his many helpful comments, and his patience
in going through various versions of this paper. Special thanks also
to {\sc Volker Kaibel} for his very careful reading, to {\sc G\"unter
  Rote} for pointing out a gap in my first version of the realization
construction, and to all members of the Discrete Geometry group at TU
Berlin for their great support and the wonderful working environment
there.

\end{document}